\documentclass[10pt]{article}

\usepackage[dvips]{graphicx}
\usepackage[T1]{fontenc}
\usepackage[utf8]{inputenc}
\usepackage{authblk}
\usepackage{color}
\usepackage{mathtools}

\usepackage[numbers]{natbib}
\usepackage[OT4]{fontenc}
\usepackage[utf8]{inputenc}

\usepackage{amsthm}
\usepackage{amsmath}
\usepackage{bbm}
\usepackage{amssymb}

\newcommand{\Exp}{{\rm I\hspace{-0.8mm}E}}
\newcommand{\Prob}{{\rm I\hspace{-0.8mm}P}}

\newcommand{\iz}{{\rm \rlap Z\kern 2.2pt Z}}

\newcommand{\pk}[1]{{\rm I\hspace{-0.8mm}P}\left(#1\right)}

  \def\td{\text{\rm d}}
\def\IF{\infty}
\newcommand{\halmos}{\vspace{3mm} \hfill $\Box$}

\newcommand{\ind}{{1\hspace{-1mm}{\rm I}}}

\newtheorem{theorem}{Theorem}

\newtheorem{corollary}{Corollary}

\newtheorem{remark}{Remark}
\newtheorem{example}{Example}

\author[]{Zbigniew Michna\footnote{zbigniew.michna@ue.wroc.pl}}
\affil[]{Department of Mathematics and Cybernetics\\
Wrocław University of Economics\\ 
53-345 Wrocław Poland}

\title{\bf Ruin probabilities for two collaborating\\ insurance companies}

\date{}

\begin{document}

\maketitle

\begin{abstract}
In this note we find a formula for the supremum distribution of spectrally positive or negative L\'evy processes
with a broken linear drift. This gives formulas for ruin probabilities in the case when two insurance companies (or two branches of the same company) divide between them both claims and premia in some specified proportions. As an example we consider gamma L\'evy process, $\alpha$-stable L\'evy process and Brownian motion. 
Moreover we obtain identities for Laplace transform of the distribution for  the supremum of L\'evy processes with randomly broken drift and on random intervals.

\vspace{5mm}
\newline
\vspace{2cm}
\end{abstract}

\section{Introduction}
In this paper we study the supremum distribution of a spectrally positive or negative L\'evy process with a piecewise linear drift. We find exact formulas for the distribution of supremum which are expressed by one-dimensional densities of a given L\'evy process. The results can be applied to find ruin probabilities in the case when two insurance companies (or two branches of the same company) divide between them both claims and premia in some specified proportions (proportional reinsurance). 
%We investigate the problem of the so-called dependent risks (see Denuit et al. \cite{de:ge:ma:99}, Dhaene and Denuit\cite{dh:de:99}, Dhaene and Goovaerts \cite{dh:go:96}, Dhaene and Goovaerts \cite{dh:go:97}). %
Moreover the formulas can be used for a two-node tandem queue (see Lieshout and Mandjes \cite{li:ma:07}). Avram et al.  \cite{av:pa:pi:08:a} investigates a spectrally positive L\'evy process with a broken drift (reduction of the risk problem to one dimension) and they find the double Laplace transform of the infinite time survival probability. As an example they obtain an explicit analytical representation of the infinite time survival probability if the claims are exponentially distributed (the compound Poisson process with exponentially distributed claims).
In Avram et al. \cite{av:pa:pi:08} the related problem is investigated if the accumulated claim amount is modeled by a L\'evy process that admits negative exponential
moments. They find exact formulas for ruin probabilities expressed by ordinary ruin probabilities when the accumulated claim amount process is spectrally negative or a compound Poisson process with exponential claims. Additionally they find asymptotic behavior of ruin probabilities under the Cram\'er assumption. In Foss et al. \cite{fo:ko:pa:ro:17} the same problem is investigated as in Avram et al. \cite{av:pa:pi:08} but the subexponential claims are admitted and an asymptotic behavior of ruin probabilities on finite and infinite time horizon is found.
In the models analyzed in this contribution we assume that the accumulated claim amount process is a spectrally positive or a spectrally negative L\'evy process with one-dimensional density functions.
%(admitting an atom at zero that is the densities can be improper I THINK WE HAVE TO REMOVE IT, SEE THE EXPLANATION BELOW)
We find exact formulas for ruin probabilities expressed  by one-dimensional densities of an underlying L\'evy process. The main difference of our models and models of Avram et al. \cite{av:pa:pi:08}  is that we admit heavy tailed claims and we provide explicit formulas of ruin probabilities both on finite and infinite time horizon unlike Avram et al.  \cite{av:pa:pi:08:a} and Avram et al. \cite{av:pa:pi:08} where it is done only on the infinite time horizon. 

The layout of the rest of the article is the following. In this section we recall the formulas which will be used in the main results. The next section contains the main results that is the distribution of  supremum of a L\'evy process with a broken drift and examples. In Section \ref{sec3} we outline how to apply the main results to ruin probabilities for two collaborating insurance companies.
The last section deals with the identities for Laplace transform of the distribution for  supremum of 
L\'evy processes with a randomly broken drift and on random intervals.

In Michna et al. \cite{mi:pa:pi:15} a joint distribution of the random variable $Y(T)$
and $\inf_{t< T}Y(t)$ was found where $Y$ is a spectrally negative L\'evy process (we will consider real stochastic processes with time defined on the non-negative half real line).
\begin{theorem}\label{tacneginf}
If $Y$ is a spectrally negative L\'evy process and the one-dimensional distributions of $Y$ are absolutely continuous then
$$
\Prob (\inf_{t< T} Y(t)< -u, \,Y(T)+u\in \td z)=
\td z\int_0^T\frac{z}{T-s}\, p(z, T-s)\,p(-u,s)\td s,
$$
where $T, u>0$, $z\geq 0$ and $p(x,s)$ is a density function of $Y(s)$ for $s>0$.
\end{theorem}
\begin{remark}
We do not expose a linear drift of the process $Y$ but it can be incorporated in the process $Y$.
\end{remark}
%\begin{remark}
%In the above theorem we can assume that the one-dimensional distributions of $Y$ can have
%an atom at zero then the density $p(x,s)$ is an improper density function. Thus
%$Y$ can be a compound Poisson process as well. $Y$ COMPOUND POISSON WITH ATOM AT ZERO HAS NO DRIFT SO IT IS NOT INTERESTING. COMPOUND POISSON WITH DRIFT HAS ATOM AT $ct$.
%\end{remark}
If $X$ is a spectrally positive L\'evy process then $X=-Y$ and we get the following corollary.
\begin{corollary}\label{jsp}
If $X$ is a spectrally positive L\'evy process and the one-dimensional distributions of $X$ are absolutely continuous then
\begin{eqnarray*}
\lefteqn{\Prob (\sup_{t< T} X(t)\leq u, \,X(T)\in \td z)=}\\
&&\left[f(z,T)-\int_0^T\frac{u-z}{T-s}\, f(z-u, T-s)\,f(u,s)\td s\right]\td z,
\end{eqnarray*}
where $T, u>0$, $z\in (-\infty, u]$ and $f(x,s)$ is a density function of $X(s)$ for $s>0$.
\end{corollary}
Integrating the last formula with respect to $z$ we get the following theorem (see Michna et al. \cite{mi:pa:pi:15} and Michna \cite{mi:11}).
\begin{theorem}\label{mi}
If the one-dimensional distributions of $X$ are absolutely continuous then
\begin{equation}\label{mi1}
\Prob(\sup_{t< T} X(t)>u)
=\Prob(X(T)>u)+\int_0^T\frac{\Exp(X(T-s))^-}{T-s}\,
f(u,s)\,{\rm d}s\,,
%\label{mainf}
\end{equation}
where $x^-=-\min\{x,0\}$ and $f(u,s)$ is a density function of $X(s)$ for $s>0$.
\end{theorem}
\begin{remark}
The above formula extends the result of Tak\'acs \cite{ta:65} to L\'evy processes with infinite variation.
\end{remark}
Let us now find the joint distribution of supremum and the value of the process for any spectrally negative L\'evy process. It will easily follow from Corollary \ref{jsp} and the duality lemma.
\begin{corollary}\label{refy}
If $Y$ is a spectrally negative L\'evy process and the one-dimensional distributions of $Y$ are absolutely continuous then
$$
\Prob (\sup_{t< T} Y(t)\leq u, \,Y(T)\in \td z)=
\left[p(z,T)-u\int_0^T\frac{p(u, T-s)}{T-s}\,p(z-u,s)\td s\right]\td z\,,
$$
where $T, u>0$, $z\in (-\infty, u]$ and $p(x,s)$ is a density function of $Y(s)$ for $s>0$.
\end{corollary}
\proof
By the duality lemma (see e.g. Bertoin \cite{be:96}) we have that\\ $X((T-t)-)-X(T)\stackrel{d}{=}Y(t)$ in the sense of finite dimensional distributions for $t\leq T$ ($X(t-)$ means the left-hand side limit at $t$). Thus we get
\begin{eqnarray*}
\lefteqn{\Prob (\sup_{t< T} X(t)\leq u, \,X(T)\in \td z)=}\\
&&\Prob (\sup_{t< T} X((T-t)-)\leq u, \,X(T)\in \td z)\\
&&=\Prob (\sup_{t< T} (X((T-t)-)-X(T))\leq u-z, \,X(T)\in \td z)\\
&&=\Prob (\sup_{t< T} Y(t)\leq u-z, \,-Y(T)\in \td z)\,.
\end{eqnarray*}
Substituting $u'=u-z$ and $z'=-z$ and using  Corollary \ref{jsp} we obtain the formula.
\halmos

Integrating the last formula with respect to $z$ we could get a similar result to eq. (\ref{mi1}) for spectrally negative L\'evy processes. However we obtain a simpler formula from Kendall's identity (see Kendall \cite{ke:57}).
The following theorem can be found in a more general form in Tak\'acs \cite{ta:65} (see also Michna \cite{mi:13} for the distribution of supremum for spectrally negative L\'evy processes).
\begin{theorem}\label{supy}
If $Y$ is a spectrally negative L\'evy process and the one-dimensional distributions of $Y$ are absolutely continuous then
$$
\Prob (\sup_{t< T} Y(t)>u)=u\int_0^T\frac{p(u,s)}{s}\,\td s\,,
$$
where $p(u,s)$ is the density function of $Y(s)$.
\end{theorem}
\proof
It follows directly from Kendall's identity (see Kendall \cite{ke:57} or e.g. Sato \cite{sa:99} Th. 46.4).
\halmos

\section{Main results and examples}
In this section we analyze the distribution of supremum for both $X(t)-c(t)$ and $Y(t)-c(t)$ where $X$ is a spectrally positive L\'evy process and $Y$ is a spectrally negative L\'evy process and
\begin{equation}\label{d}
c(t)=
\left\{\begin{array}{ll}
c_1 t &\mbox{if }\, t\in[0, T]\\
c_2(t-T)+c_1 T&\mbox{if }\, t\in(T, \infty)\,,
\end{array}
\right.
\end{equation}
where $c_1, c_2\geq 0$. %and $X(t)-c_i t$, $Y(t)-c_i t$ are not monotone for $i=1,2$.
Since we now expose the drift of the process we will assume that densities of $X(s)$ and $Y(s)$ are
$f(x,s)$ and $p(x,s)$, respectively (unlike the previous section where a linear drift could be incorporated in the processes). 
%Let us also notice that the considered L\'evy processes have now a linear drift (e.g. $X$ can be a subordinator).
\begin{theorem}\label{main}
If $S>T$ ($S$ is finite or $S=\infty$) and $X(t)$ is absolutely continuous with density $f(x,t)$ then
\begin{eqnarray*}
\lefteqn{\Prob(\sup_{t<S}(X(t)-c(t))>u)=A+B\coloneqq}\\
&&\Prob(\sup_{t< T}(X(t)-c_1 t)>u)\\
&&+\,\Prob(\sup_{t< T}(X(t)-c_1 t)\leq u, \sup_{0<t<S-T}(X(t+T)-X(T)-c_2 t)>u-X(T)+c_1 T)\,,
\end{eqnarray*}
where
$$
A=\Prob(X(T)-c_1 T>u)+\int_0^T\frac{\Exp(X(T-s)-c_1 (T-s))^-}{T-s}\,
f(u+c_1 s,s)\,{\rm d}s
$$
and
\begin{eqnarray*}
B&=&\int_{0}^\infty \Prob(\sup_{t< S-T}(X(t)-c_2 t)>z)f(-z+u+c_1 T,T)\td z\\
&&-\int_{0}^\infty z\,\Prob(\sup_{t< S-T}(X(t)-c_2 t)>z)\td z\\
&&\,\,\,\,\,\,\cdot\int_{0}^T\frac{f(u+c_1 s,s)}{T-s}f(-z+c_1(T-s),T-s)\td s\,.
\end{eqnarray*}
\end{theorem}
\proof
The decomposition $A+B$ we get as follows
\begin{eqnarray*}
\lefteqn{\Prob(\sup_{t<S}(X(t)-c(t))>u)=A+B\coloneqq}\\
&&\Prob(\sup_{t< T}(X(t)-c_1 t)>u)\\
&&+\,\Prob((\sup_{t< T}(X(t)-c_1 t)\leq u, \sup_{T<t<S}(X(t)-c_2(t-T)-c_1T)>u)\\
&&=\Prob(\sup_{t< T}(X(t)-c_1 t)>u)\\
&&+\,\Prob(\sup_{t< T}(X(t)-c_1 t)\leq u, \sup_{0<t<S-T}(X(t+T)-X(T)-c_2 t)>u-X(T)+c_1 T)\,.
\end{eqnarray*}
The formula for $A$ we directly get from Theorem \ref{mi}.
Let $F(\td x, \td z)$ be the joint distribution of $(\sup_{t< T}(X(t)-c_1 t), X(T)-c_1 T)$.
Then the formula for $B$ follows from the strong Markov property and Corollary \ref{jsp} that is
\begin{eqnarray*}
B&=&\int_0^u\int_{-\infty}^u\Prob(\sup_{t<S-T}(X(t)-c_2 t)>u-z)F(\td x, \td z)\\
&=&\int_{-\infty}^u\Prob(\sup_{t<S-T}(X(t)-c_2 t)>u-z)f(z+c_1 T, T)\td z\\
&&-\int_{-\infty}^u\Prob(\sup_{t<S-T}(X(t)-c_2 t)>u-z)\td z\\
&&\cdot\int_0^T\frac{u-z}{T-s}f(z-u+c_1 (T-s), T-s)f(u+c_1s, s)\td s
\end{eqnarray*}
and substituting $z'=u-z$ we obtain the final formula.
\halmos

Similarly we get a formula for spectrally negative L\'evy processes.
\begin{theorem}
If $S>T$ ($S$ is finite or $S=\infty$) and $Y(t)$ is absolutely continuous with density $p(x,t)$ then
$\Prob(\sup_{t<S}(Y(t)-c(t))>u)=A+B$
where
$$
A=\Prob(\sup_{t< T}(Y(t)-c_1 t)>u)=u\int_0^T\frac{p(u+c_1s,s)}{s}\,\td s
$$
and
\begin{eqnarray*}
B&=&\int_{0}^\infty \Prob(\sup_{t< S-T}(Y(t)-c_2 t)>z)p(-z+u+c_1 T,T)\td z\\
&&-u\,\int_{0}^\infty \Prob(\sup_{t< S-T}(Y(t)-c_2 t)>z)\td z\\
&&\,\,\,\,\,\,\cdot\int_{0}^T\frac{p(-z+c_1s, s)}{T-s}\,p(u+c_1(T-s), T-s)\td s\,.
\end{eqnarray*}
\end{theorem}
\proof
Using Corollary \ref{refy} and Th. \ref{supy} we proceed the same way as in the proof of Th. \ref{main}.
\halmos

The application of Th. \ref{main} leads to the following example with Brownian motion (see Mandjes \cite{ma:04} and Lieshout and Mandjes \cite{li:ma:07} or Avram et al. \cite{av:pa:pi:08}).
\begin{example}\label{brexpl}
If $W$ is the standard Brownian motion then
\begin{eqnarray*}
\lefteqn{\Prob(\sup_{t<\infty}(W(t)-c(t))>u)=}\\
&&\Phi(-uT^{-1/2}-c_1\sqrt{T})+e^{-2c_1u}\Phi(-uT^{-1/2}+c_1\sqrt{T})\\
&&\,\,\,\,\,\,+e^{-2c_2(u+c_1T-c_2T)}\Phi(uT^{-1/2}+(c_1-2c_2)\sqrt{T})\\
&&\,\,\,\,\,\,-e^{2(c_2-c_1)u+2c_2^2T-2c_1c_2T}\Phi(-uT^{-1/2}+(c_1-2c_2)\sqrt{T})\,.
\end{eqnarray*}
%\end{corollary}
%\proof
Indeed using Theorem \ref{main} and
$$
\Prob(\sup_{t< T}(W(t)-ct)>u)=\Phi(-uT^{-1/2}-c\sqrt{T})+e^{-2cu}\Phi(-uT^{-1/2}+c\sqrt{T})
$$
and
$$
\Prob(\sup_{t<\infty}(W(t)-ct)>u)=e^{-2cu}
$$
for $c\geq 0$ (see e.g. Dębicki and Mandjes \cite{de:ma:15})  we get
\begin{equation}\label{abbr}
\Prob(\sup_{t<\infty}(W(t)-c(t))>u)=A+B\,,
\end{equation}
where
\begin{equation}\label{ab}
A=A(c_1, T, u)\coloneqq\Phi(-uT^{-1/2}-c_1\sqrt{T})+e^{-2uc_1}\Phi(-uT^{-1/2}+c_1\sqrt{T})
\end{equation}
and
\begin{eqnarray*}
\lefteqn{B=}\\
&&e^{-2c_2(u+c_1T-c_2T)}\Phi(uT^{-1/2}+(c_1-2c_2)\sqrt{T})\\
&&\,\,\,\,\,-\,\frac{e^{-c_1u-c_1^2T/2}}{2\pi}\int_{0}^\infty ze^{(c_1-2c_2)z}\td z\int_0^T(T-s)^{-3/2}s^{-1/2}e^{-\frac{z^2}{2(T-s)}-\frac{u^2}{2s}}\td s\,.
\end{eqnarray*}
Let us take $c=c_1=c_2\geq 0$ in eq. (\ref{abbr}). Then $A+B=e^{-2uc}$ and the second summand
of $A$ and the first one of $B$ sum up to $e^{-2uc}$ thus we get
$$
\frac{e^{-cu-c^2T/2}}{2\pi}\int_{0}^\infty ze^{-cz}\td z\int_0^T(T-s)^{-3/2}s^{-1/2}e^{-\frac{z^2}{2(T-s)}-\frac{u^2}{2s}}\td s=
\Phi(-uT^{-1/2}-c\sqrt{T})\,.
$$
Thus using the last identity for $c=2c_2-c_1\geq 0$ we get the second therm of $B$.

Similarly let us take $c=c_1$ and $c_2=0$ in eq. (\ref{abbr}). Then $A+B=1$ and the first summand
of $A$ and the first one of $B$ sum up to 1 thus we get
$$
\frac{e^{cu-c^2T/2}}{2\pi}\int_{0}^\infty ze^{cz}\td z\int_0^T(T-s)^{-3/2}s^{-1/2}e^{-\frac{z^2}{2(T-s)}-\frac{u^2}{2s}}\td s=
\Phi(-uT^{-1/2}+c\sqrt{T})\,.
$$
Thus using the last identity for $c=c_1-2c_2> 0$ we get the second therm of $B$.
%\halmos
\end{example}
\begin{example}
Let $0<T<S<\infty$ and $W$ be the standard Brownian motion then
\begin{eqnarray*}
\lefteqn{\Prob(\sup_{t<S}(W(t)-c(t))>u)=}\\
&&A(c_1, T, u)
+\frac{1}{\sqrt{2\pi T}}\int_0^\infty A(c_2, S-T, z) e^{-\frac{(u+c_1T-z)^2}{2T}}\td z\\
&&-\frac{e^{-uc_1-\frac{c_1^2 T}{2}}}{2\pi}\int_0^\infty ze^{c_1 z} A(c_2, S-T, z)\td z
\int_0^T s^{-1/2}(T-s)^{-3/2}\,e^{-\frac{z^2}{2(T-s)}-\frac{u^2}{2s}}\td s\,,
\end{eqnarray*}
where $A(c_1, T, u)$ is defined in eq. (\ref{ab}).
\end{example}
\begin{example}
Let $X(t)$ be gamma L\'evy process with the density
$$f(x, t)=\frac{\delta^t}{\Gamma(t)}x^{t-1}e^{-\delta x}\ind_{\{x>0\}}$$
where $\delta>0$ and $c(t)$ be defined in eq. (\ref{d}). Using Th. \ref{main} we give the explicit formulas of $\Prob(\sup_{t<S}(X(t)-c(t))>u)=A+B$  for both $T<S<\IF$ and $S=\IF$, respectively. For $T<S<\IF$, we have that
\begin{eqnarray*}
A&=&\frac{\delta^T}{\Gamma(T)}\int_{u+c_1T}^\IF x^{T-1}e^{-\delta x}\td x\\
&&+\,\delta^T e^{-\delta u}\int_0^T\frac{(u+c_1s)^{s-1}e^{-c_1\delta s }}{\Gamma(s)\Gamma(T-s+1)}\td s\int_0^{c_1(T-s)}(c_1(T-s)-x)x^{T-s-1}e^{-\delta x} \td x\\
&\eqqcolon&A(c_1,T,u)
\end{eqnarray*}
and
\begin{eqnarray*}
\lefteqn{B=}\\
&&\frac{\delta^S e^{-\delta (u+c_1T)}}{\Gamma(T)\Gamma(S-T)}\int_0^{u+c_1T} (u+c_1T-z)^{T-1}e^{\delta z}\td z\int_{z+c_2(S-T)}^\IF x^{S-T-1}e^{-\delta x}\td x\\
&&+\frac{\delta^S e^{-\delta (u+c_1T)}}{\Gamma(T)}\int_0^{u+c_1T} (u+c_1T-z)^{T-1}\td z\int_0^{S-T}\frac{(z+c_2s)^{s-1}e^{-c_2\delta s}}{\Gamma(s)\Gamma(S-T-s+1)}\td s\\
&&\quad \cdot \int_0^{c_2(S-T-s)}(c_2(S-T-s)-x)x^{S-T-s-1}e^{-\delta x}\td x\\
&&-\delta^T e^{-\delta (u+c_1T)}\int_0^{c_1T} z e^{\delta z}A(c_2,S-T,z)\td z\\
&&\cdot\int_0^{\frac{c_1T-z}{c_1}}\frac{(u+c_1s)^{s-1}(c_1(T-s)-z)^{T-s-1}}{\Gamma(s)\Gamma(T-s+1)}\td s\,.
\end{eqnarray*}
For $S=\IF$, we additionally assume that $c_2\delta >1$. In this case, since $X(t)$ has finite variation, in view of Th. 4 in Tak\'acs \cite{ta:65} we have
$$\Prob(\sup_{t<\IF}(X(t)-c_2t)>z)=\frac{c_2\delta-1}{\delta}\,e^{-\delta z}\int_0^\IF\frac{\delta^s}{\Gamma(s)}(z+c_2 s)^{s-1}e^{-\delta c_2 s}\td s\,, \,z>0\,.$$
Let us notice that $A$ is the same as in the case $T<S<\IF$ and using the above expression we get
\begin{eqnarray*}
\lefteqn{B=}\\
&&\frac{(c_2\delta-1)\delta^{T-1}e^{-\delta(u+c_1T)}}{\Gamma(T)}\int_0^{u+c_1T}(u+c_1T-z)^{T-1}\td z\\ 
&&\cdot\int_0^\IF\frac{\delta^s}{\Gamma(s)}(z+c_2s)^{s-1}e^{-\delta c_2 s}\td s\\
&&-\,(c_2\delta-1)\delta^{T-1} e^{-\delta (u+c_1T)}\int_0^{c_1T} z \td z\\
&&\cdot\int_0^{T-\frac{z}{c_1}}\frac{(u+c_1s)^{s-1}(c_1(T-s)-z)^{T-s-1}}{\Gamma(s)\Gamma(T-s+1)}\td s\int_0^\IF\frac{\delta^t}{\Gamma(t)}(z+c_2 t)^{t-1}e^{-\delta c_2 t}\td t.
\end{eqnarray*}
\end{example}
\begin{example}
Let $Z(s)$ be an $\alpha$-stable L\'evy process totally skewed to the right (that is with $\beta=1$ see e.g. Janicki and Weron \cite{ja:we:94} or Samorodnitsky and Taqqu \cite{sa:ta:94}) with $1<\alpha<2$ and expectation zero then its density function is the following
$$
f(x,s)=\frac{1}{\pi s^{1/\alpha}}
\int_0^\infty e^{-t^\alpha}\cos\left(ts^{-1/\alpha}x-t^\alpha\tan{\frac{\pi\alpha}{2}}\right)\td t
$$
(see e.g. Nolan \cite{no:97}). Then (see Michna et al. \cite{mi:pa:pi:15}) for $c>0$
\begin{eqnarray*}
\lefteqn{A(c,\infty,u)\coloneqq\Prob(\sup_{t<\infty}(Z(t)-ct)>u)=}\\
&&\frac{c}{\pi} \int_0^\infty s^{-1/\alpha}\,\td s
\int_0^\infty e^{-t^\alpha}\cos\left(ts^{-1/\alpha}(u+cs)-t^\alpha\tan{\frac{\pi\alpha}{2}}\right)\td t
\end{eqnarray*} and for any $c$ and $T>0$
\begin{eqnarray}
\lefteqn{A(c,T,u)\coloneqq \Prob(\sup_{t<T}(Z(t)-ct)>u)=}\label{astable}\\
&&\frac{1}{\pi T^{1/\alpha}} \int_u^\infty\td x
\int_0^\infty e^{-t^\alpha}\cos\left(tT^{-1/\alpha}(x+cT)-t^\alpha\tan{\frac{\pi\alpha}{2}}\right)\td t\nonumber\\
&& +\frac{1}{\pi}\int_0^T \frac{\Exp(Z(T-s)-c(T-s))^-}{(T-s)s^{1/\alpha}}\,\td s\nonumber\\
&&\,\,\,\,\,\,\,\,\cdot\int_0^\infty e^{-t^\alpha}\cos\left(ts^{-1/\alpha}(u+cs)-t^\alpha\tan{\frac{\pi\alpha}{2}}\right)\td t\nonumber
\end{eqnarray}
where
$$
\Exp(Z(s)-cs)^-=\frac{-1}{\pi s^{1/\alpha}}\int_{-\infty}^0 x\, \td x\int_0^\infty e^{-t^\alpha}\cos\left(ts^{-1/\alpha}(x+cs)-t^\alpha\tan{\frac{\pi\alpha}{2}}\right)\td t\,.
$$
Thus using Th. \ref{main} for $S>T>0$ (allowing also $S=\infty$ and putting $\infty-T=\infty$) we get
$$
\Prob(\sup_{t<S}(Z(t)-c(t))>u)=A+B=A+B_1-B_2\,,
$$
where $A=A(c_1, T, u)$ (see eq. (\ref{astable})) and
\begin{eqnarray*}
\lefteqn{B_1=}\\
&&\frac{1}{\pi T^{1/\alpha}}\int_0^\infty A(c_2,S-T,z)\,\td z\\
&&\cdot\int_0^\infty e^{-t^\alpha}\cos\left(tT^{-1/\alpha}(-z+u+c_1T)-t^\alpha\tan{\frac{\pi\alpha}{2}}\right)\td t
\end{eqnarray*}
and
\begin{eqnarray*}
\lefteqn{B_2=}\\
&&\frac{1}{\pi^2}\int_0^\infty z\, A(c_2, S-T, z)\,\td z\int_0^T\frac{\td s}{(T-s)^{1/\alpha+1}s^{1/\alpha}}\\
&&\cdot\int_0^\infty  e^{-t^\alpha}\cos\left(ts^{-1/\alpha}(u+c_1s)-t^\alpha\tan{\frac{\pi\alpha}{2}}\right)\td t\\
&& \cdot\int_0^\infty e^{-w^\alpha}\cos\left(w(T-s)^{-1/\alpha}(-z+c_1(T-s))-w^\alpha\tan{\frac{\pi\alpha}{2}}\right)\td w\,.
\end{eqnarray*}
\end{example}

\section{Two collaborating insurance companies}\label{sec3}
Let us consider two insurance companies which split the amount they pay out of each claim in
proportions $\delta_1> 0$ and $\delta_2> 0$ where $\delta_1+\delta_2=1$, and receive premiums at rates $p_1>0$ and $p_2>0$, respectively (see Avram et al. \cite{av:pa:pi:08}). Then the corresponding risk processes are
$$
R_i(t)=x_i+p_i t-\delta_i X(t)\,,
$$
where $i=1,2$, $x_i>0$ and $X(t)$ is an accumulated claim amount up to time $t$. One can be interested in the instant when at least one insurance company is ruined
$$
\tau_{or}(x_1, x_2)=\inf\{t>0: R_1(t)<0\,\, \mbox{or}\,\, R_2(t)<0\}
$$
and in the instant when both insurance companies are simultaneously ruined
$$
\tau_{sim}(x_1, x_2)=\inf\{t>0: R_1(t)<0\,\, \mbox{and}\,\, R_2(t)<0\}\,.
$$
Let the ultimate ruin probabilities be
$$
\psi_{or}(x_1, x_2)=\Prob(\tau_{or}(x_1, x_2)<\infty)\,,\,\,\,\,\,\,
\psi_{sim}(x_1, x_2)=\Prob(\tau_{sim}(x_1, x_2)<\infty)
$$
and
$$
\psi_{1}(x_1)=\Prob(\tau_1(x_1)<\infty)\,,\,\,\,\,\,\,
\psi_{2}(x_2)=\Prob(\tau_2(x_2)<\infty)\,,
$$
where $\tau_i(x_i)=\inf\{t>0: R_i(t)<0\}$ for $i=1,2$. One can also be interested in the following
ruin probability
$$
\psi_{and}(x_1, x_2)=\Prob(\tau_1(x_1)<\infty\,\,\mbox{and}\,\,\tau_2(x_2)<\infty)
$$
and the following relation is easily to notice
$$
\psi_{and}(x_1, x_2)=\psi_1(x_1)+\psi_2(x_2)-\psi_{or}(x_1,x_2)\,.
$$
Let us put $u_i=x_i/\delta_i$ and $c_i=p_i/\delta_i$ where $i=1,2$. Then we get
$$
\tau_{or}(x_1, x_2)=\inf\{t>0: X(t)>u_1+c_1 t\,\,\mbox{or}\,\,X(t)>u_2+c_2 t\}
$$
and
$$
\tau_{sim}(x_1, x_2)=\inf\{t>0: X(t)>u_1+c_1 t\,\,\mbox{and}\,\,X(t)>u_2+c_2 t\}\,.
$$
If the lines $y=u_1+c_1 t$ and $y=u_2+c_2 t$ do not cross each other in the first quadrant  then
the ruin probabilities $\psi_{or}(x_1, x_2)$ and $\psi_{sim}(x_1, x_2)$ reduce to ordinary ruin probabilities of a risk process with a linear drift. If they cross each other in the first quadrant and e.g. $u_1<u_2$ ($c_1>c_2$) then
\begin{equation}\label{orfor}
\psi_{or}(x_1, x_2)=\Prob(\sup_{t<\infty}(X(t)-c(t))>u_1)\,,
\end{equation}
where $c(t)$ is defined in eq. (\ref{d}) with $T=(u_2-u_1)/(c_1-c_2)$ (we take $c(t)=\min(u_1+c_1 t, u_2+c_2 t)-u_1$).

Similarly, if the lines have a common point in the first quadrant and  e.g. $u_2<u_1$ ($c_2>c_1$) then
\begin{equation}\label{simfor}
\psi_{sim}(x_1, x_2)=\Prob(\sup_{t<\infty}(X(t)-c(t))>u_1)\,,
\end{equation}
where $c(t)$ is defined in  eq. (\ref{d}) with $T=(u_2-u_1)/(c_1-c_2)$ (we take $c(t)=\max(u_1+c_1 t, u_2+c_2 t)-u_1$).

\begin{example}
Let $X(t)$ be the standard Brownian motion. Then using eq. (\ref{orfor}) and Example \ref{brexpl}
we get for $u_1<u_2$ and $c_1>c_2$
\begin{eqnarray*}
\lefteqn{\psi_{or}(x_1, x_2)=}\\
&&\Phi(a(-u_1,-c_1))+e^{-2c_1u_1}\Phi(a(-u_1, c_1))\\
&&\,\,\,+e^{-2c_2u_2}\Phi(a(u_1,  c_1-2c_2))-e^{-2(c_1-2c_2)u_1-2c_2u_2}\Phi(a(-u_1, c_1-2c_2))\,,
\end{eqnarray*}
where $a(u,c)=uT^{-1/2}+c\sqrt{T}$, $T=(u_2-u_1)/(c_1-c_2)$, $u_i=x_i/\delta_i$ and $c_i=p_i/\delta_i$ for $i=1,2$. This formula recovers the result of Avram et al. \cite{av:pa:pi:08} Eq. (67).
The same way we obtain the formula for $\psi_{sim}(x_1, x_2)$.
\end{example}
In a similar way we can consider ruin probabilities on a finite time horizon.

\section{Randomly broken drift and random\\ interval}
In the fluctuation theory there are many interesting identities for L\'evy processes and an exponentially distributed time e.g. the distribution of supremum on an exponentially distributed time interval (see e.g. Bertoin \cite{be:96} Sec. VI. 2. and Sec. VII).
Thus let us consider a spectrally positive L\'evy process $X$ with a randomly broken drift that is let us assume that $T$ (see eq. (\ref{d})) is an exponentially distributed random variable with mean $1/\lambda$ independent of the process $X$.
Moreover let us investigate two cases $S=\infty$ (see Th. \ref{main}) and $S-T=V$ is a positive random variable independent of the process $X$ and the random variable
$T$.
We put
$$\varphi_i(\gamma)=\ln \Exp\exp(-\gamma( X(1)-c_i)),\, i=1,2\,,$$
where $\gamma\geq 0$
and $\overleftarrow{\varphi}_i(\lambda), i=1,2$ is the inverse function of $\varphi_i$.
\begin{theorem}\label{mainlap}
If $X$ is a spectrally positive L\'{e}vy process and $T$ is an exponential random variable with mean $1/\lambda>0$ independent of $X$ then for any  \\
$\gamma>\overleftarrow{\varphi}_1(\lambda)$
\begin{eqnarray}\label{mainexp}
\lefteqn{\Exp e^{-\gamma\sup_{t< T+V} (X(t)-c(t))}=}\nonumber\\
&&\Exp e^{-\gamma\sup_ {t< T} (X(t)-c_1t)}\\
&&+\,\frac{\gamma\lambda}{\varphi_1(\gamma)-\lambda}\left[\frac{1-\Exp e^{-\gamma\sup_{ t< V}(X(t)-c_2 t)}}{\gamma}-\frac{1-\Exp e^{-\overleftarrow{\varphi}_1(\lambda)\sup_{t<V}(X(t)-c_2 t)}}{\overleftarrow{\varphi}_1(\lambda)}\right]\nonumber
\end{eqnarray}
where $V$ is a positive random variable  independent of $X$ and $T$.
\end{theorem}
\proof
Observe that for $\gamma>0$
 \begin{equation*}\Exp e^{-\gamma\sup_{t<T+V} (X(t)-c(t))}=1-\gamma\int_{0}^\IF e^{-\gamma u}\pk{\sup_{t<T+V}X(t)-c(t)>u}\td u
\end{equation*}
and
\begin{eqnarray*}
\lefteqn{\pk{\sup_{t<T+V}(X(t)-c(t))>u}}\\
&&=\pk{\sup_{t<T}(X(t)-c_1 t)>u}\\
&&+\,\pk{\sup_{t<T}(X(t)-c_1 t)\leq u, \sup_{t<V}(X(t+T)-X(T)-c_2 t)>u-X(T)+c_1 T}.
\end{eqnarray*}
Thus we have that
\begin{eqnarray}\label{I123}
\int_{0}^\IF e^{-\gamma u}\pk{\sup_{t<T+V}X(t)-c(t)>u}\td u= I_1+I_2\,,
\end{eqnarray}
where
$$I_1\coloneqq\int_{0}^\IF e^{-\gamma u}\pk{\sup_{t<T}X(t)-c_1t>u}\td u=\frac{1-\Exp e^{-\gamma\sup_{t<T} (X(t)-c_1t)}}{\gamma}$$
and
\begin{eqnarray*}
\lefteqn{I_2 \coloneqq\int_{0}^\IF e^{-\gamma u}}\\
&&\cdot\pk{\sup_{t<T}(X(t)-c_1 t)\leq u, \sup_{t<V}(X(t+T)-X(T)-c_2 t)>u-X(T)+c_1 T}\td u\,.
\end{eqnarray*}
By the fact that  $T$ is exponentially distributed and independent of $X$ and $V$ we have
\begin{eqnarray*}
\lefteqn{I_2=}\\
&&\lambda\int_0^\IF e^{-\lambda s}\td s\,\int_{0}^\IF e^{-\gamma u}\\
&&\cdot\pk{\sup_{t<s}(X(t)-c_1 t)\leq u, \sup_{t<V}(X(t+s)-X(s)-c_2 t)>u-X(s)+c_1 s}\td u\,.
\end{eqnarray*}
Moreover, by the independence of $X(t+s)-X(s)-c_2 t,\, t\geq 0$ and $X(s)-c_1s$ and the fact that
$$\pk{\sup_{t<s}(X(t)-c_1 t)\leq u, u-X(s)+c_1 s\leq z}=0\,, \quad z<0$$
we have
\begin{eqnarray*}
I_2&=&\lambda\int_0^\IF e^{-\lambda s}\,\td s\int_{0}^\IF e^{-\gamma u}\,\td u\int_{0}^\IF \pk{\sup_{t<V}(X(t)-c_2 t)>z}\\
&&\cdot\pk{\sup_{t<s}(X(t)-c_1 t)\leq u, u-X(s)+c_1 s\in \td z}\\
&=&\lambda\int_{0}^\IF e^{-\gamma u}\,\td u\int_{0}^\IF e^{-\lambda s}\,\td s\int_0^\IF\pk{\sup_{t<V}(X(t)-c_2 t)>z}\\
&&\cdot \pk{\inf_{t<s}(u-X(t)+c_1 t)>0, u-X(s)+c_1 s\in \td z}\,.
\end{eqnarray*}
 Due to Suprun \cite{su:76} (see also Bertoin \cite{be:97} Lemma 1) we have that
 \begin{eqnarray*}
 &&\int_0^\IF e^{-\lambda s}\pk{\inf_{t<s}(u-X(t)+c_1 t)>0, u-X(s)+c_1 s\in \td z} \td s\\
 &&\quad = \left[e^{-\overleftarrow{\varphi}_1(\lambda)z} W^{(\lambda)}(u)-\ind(u\geq z)W^{(\lambda)}(u-z)\right]\td z\,,
 \end{eqnarray*}
where $\ind(\cdot)$ is the indicator function and $W^{(\lambda)}: [0,\IF)\rightarrow [0,\IF)$ is a continuous and increasing function such that
 $$\int_0^\IF e^{-\gamma x} \,W^{(\lambda)}(x)\td x=\frac{1}{\varphi_1(\gamma)-\lambda}\,, \quad \gamma>\overleftarrow{\varphi}_1(\lambda)\,. $$
 Consequently, for $\gamma>\overleftarrow{\varphi}_1(\lambda)$
\begin{eqnarray*}
\lefteqn{I_2=}\\
&&\lambda\int_{0}^\IF\int_{0}^\IF  e^{-\gamma u}\pk{\sup_{t<V}(X(t)-c_2 t)>z}\\
&&\cdot\left[e^{-\overleftarrow{\varphi}_1(\lambda)z} W^{(\lambda)}(u)-\mathbb{I}(u\geq z)W^{(\lambda)}(u-z)\right]\td z\td u\\
&=&\lambda \int_{0}^\IF  e^{-\overleftarrow{\varphi}_1(\lambda)z}\pk{\sup_{t< V}(X(t)-c_2 t)>z}\td z\int_{0}^\IF e^{-\gamma u}W^{(\lambda)}(u)\td u\\
&&-\lambda \int_{0}^\IF  \pk{\sup_{t<V}(X(t)-c_2 t)>z}\td z\int_0^\IF \mathbb{I}(u\geq z)e^{-\gamma u}W^{(\lambda)}(u-z)\td u \\
&=& \frac{\lambda}{\varphi_1(\gamma)-\lambda}\left[\int_{0}^\IF  e^{-\overleftarrow{\varphi}_1(\lambda)z}\pk{\sup_{t<V}(X(t)-c_2 t)>z}\td z\right.\\
&&\left.- \int_{0}^\IF  e^{-\gamma z}\pk{\sup_{t<V}(X(t)-c_2 t)>z}\td z\right]\\
&=& \frac{\lambda}{\varphi_1(\gamma)-\lambda}\left[\frac{1-\Exp e^{-\overleftarrow{\varphi}_1(\lambda)\sup_{t<V}(X(t)-c_2 t)}}{\overleftarrow{\varphi}_1(\lambda)}-\frac{1-\Exp e^{-\gamma\sup_{t<V}(X(t)-c_2 t)}}{\gamma}\right].
\end{eqnarray*}
\halmos
\begin{corollary} Under the assumption of Theorem \ref{mainlap}, if $V=\IF$, then
\begin{equation}\label{laplaceinf}
\Exp e^{-\gamma\sup_{t<\IF} (X(t)-c(t))}
=\frac{\gamma\lambda \varphi_2'(0)[\varphi_2(\gamma)-\varphi_2(\overleftarrow{\varphi}_1(\lambda))]}{\varphi_2(\gamma)(\varphi_1(\gamma)-\lambda)\varphi_2(\overleftarrow{\varphi}_1(\lambda))}\,.
\end{equation}
If $V$ is an exponentially distributed random variable with mean $1/\theta>0$ independent of $X$ and $T$ then
\begin{equation}\label{laplaceT}
\Exp e^{-\gamma\sup_{t<T+V} (X(t)-c(t))}
=
\gamma \lambda\theta\,\frac{\frac{\overleftarrow{\varphi}_2(\theta)-\overleftarrow{\varphi}_1(\lambda)}{\theta-\varphi_2(\overleftarrow{\varphi}_1(\lambda))}
-\frac{\overleftarrow{\varphi}_1(\lambda)[\overleftarrow{\varphi}_2(\theta)-\gamma]}{\gamma[\theta-\varphi_2(\gamma)]}}{\overleftarrow{\varphi}_1(\lambda)\overleftarrow{\varphi}_2(\theta)[\varphi_1(\gamma)-\lambda]}\,.
\end{equation}
\end{corollary}
\proof
{\underline{Case $V=\IF$}}. It is well-known that
\begin{eqnarray}\label{laplace1}\Exp e^{-\gamma\sup_{t<T} (X(t)-c_1t)}=\frac{\lambda}{\lambda-\varphi_1(\gamma)}\left(1-\frac{\gamma}{\overleftarrow{\varphi}_1(\lambda)}\right)\,,
\end{eqnarray}
where $\gamma\geq 0$
(see e.g. Bertoin \cite{be:96} eq. (3) p. 192 or Th. 4.1 in Dębicki and Mandjes \cite{de:ma:15}).
Moreover, by Th. 3.2 in  Dębicki and Mandjes \cite{de:ma:15} (or going with $\lambda$ to 0 in the previous identity), it follows that
$$\Exp \exp\left(-\gamma\sup_{t<\IF}(X(t)-c_2 t)\right)=\frac{\gamma\varphi_2'(0)}{\varphi_2(\gamma)}\,.$$
Consequently, by (\ref{mainexp}) for $\gamma>0$
\begin{eqnarray*}
\lefteqn{\Exp e^{-\gamma\sup_{t<\IF} (X(t)-c(t))}=}\\
&&\frac{\lambda}{\lambda-\varphi_1(\gamma)}\left[1-\frac{\gamma}{\overleftarrow{\varphi}_1(\lambda)}\right]
+\frac{\gamma\lambda}{\varphi_1(\gamma)-\lambda}\left[\frac{1-\frac{\gamma\varphi_2'(0)}{\varphi_2(\gamma)}}
{\gamma}-\frac{1-\frac{\overleftarrow{\varphi}_1(\lambda)\varphi_2'(0)}{\varphi_2(\overleftarrow{\varphi}_1(\lambda))}}{\overleftarrow{\varphi}_1(\lambda)}\right]\\
&=&\frac{\gamma \lambda\varphi_2'(0)[\varphi_2(\gamma)-\varphi_2(\overleftarrow{\varphi}_1(\lambda))]}{\varphi_2(\gamma)(\varphi_1(\gamma)-\lambda)\varphi_2(\overleftarrow{\varphi}_1(\lambda))}\,.
\end{eqnarray*}
{\underline{Case $V$ exponentially distributed}}.
Using (\ref{laplace1}), for $\gamma\geq 0$ we have that
$$\Exp\exp\left(-\gamma\sup_{t< V}(X(t)-c_2 t)\right)=\frac{\theta}{\theta-\varphi_2(\gamma)}\left(1-\frac{\gamma}{\overleftarrow{\varphi}_2(\theta)}\right).$$
Recalling (\ref{laplace1}), for $\gamma>0$ it follows that 
\begin{eqnarray*}
\lefteqn{\Exp e^{-\gamma\sup_{t<T+V} X(t)-c(t)}=}\\
&&\frac{\lambda}{\lambda-\varphi_1(\gamma)}\left(1-\frac{\gamma}{\overleftarrow{\varphi}_1(\lambda)}\right)\\
&&+\frac{\gamma\lambda}{\varphi_1(\gamma)-\lambda}\left[\frac{1-\frac{\theta}{\theta-\varphi_2(\gamma)}\left[1-\frac{\gamma}{\overleftarrow{\varphi}_2(\theta)}\right]}
{\gamma}-\frac{1-\frac{\theta}{\theta-\varphi_2(\overleftarrow{\varphi}_1(\lambda))}\left[1-\frac{\overleftarrow{\varphi}_1(\lambda)}{\overleftarrow{\varphi}_2(\theta)}\right]}
{\overleftarrow{\varphi}_1(\lambda)}\right]\\
&=&
\gamma \lambda\theta\,\frac{\frac{\overleftarrow{\varphi}_2(\theta)-\overleftarrow{\varphi}_1(\lambda)}{\theta-\varphi_2(\overleftarrow{\varphi}_1(\lambda))}
-\frac{\overleftarrow{\varphi}_1(\lambda)[\overleftarrow{\varphi}_2(\theta)-\gamma]}{\gamma[\theta-\varphi_2(\gamma)]}}{\overleftarrow{\varphi}_1(\lambda)\overleftarrow{\varphi}_2(\theta)[\varphi_1(\gamma)-\lambda]}\,.
\end{eqnarray*}
\halmos

\begin{corollary}
Let $W$ be the standard Brownian motion. Then $$\varphi_1(\gamma)=\frac{1}{2}\gamma^2+c_1\gamma\,, \quad \varphi_2(\gamma)=\frac{1}{2}\gamma^2+c_2\gamma\,,$$  
$$\overleftarrow{\varphi}_1(\lambda)=\sqrt{c_1^2+2\lambda}-c_1\,, \quad \overleftarrow{\varphi}_2(\lambda)=\sqrt{c_2^2+2\lambda}-c_2\,.$$
Consequently, for $\gamma>\sqrt{c_1^2+2\lambda}-c_1$
\begin{eqnarray*}
\lefteqn{\Exp e^{-\gamma\sup_{t<\IF} (W(t)-c(t))}=}\\
&&
\frac{\gamma \lambda c_2 (\frac{1}{2}\gamma^2+c_2\gamma-c_1^2-\lambda -(c_2-c_1)\sqrt{c_1^2+2\lambda}+c_1c_2)}{(\frac{1}{2}\gamma^2+c_1\gamma-\lambda)(\frac{1}{2}\gamma^2+c_2\gamma)(c_1^2+\lambda +(c_2-c_1)\sqrt{c_1^2+2\lambda}-c_1c_2)}\\
\end{eqnarray*}
and
\begin{eqnarray*}
\lefteqn{\Exp e^{-\gamma\sup_{t<T+V} (W(t)-c(t))}=}\\
&&
\gamma\lambda\theta\, \frac{\frac{\sqrt{c_2^2+2\theta}-\sqrt{c_1^2+2\lambda}+c_1-c_2}{\theta-c_1^2-\lambda-(c_2-c_1)\sqrt{c_1^2+2\lambda}+c_1c_2}-\frac{(\sqrt{c_1^2+2\lambda}-c_1)(\sqrt{c_2^2+2\theta}-c_2-\gamma)}{\gamma(\theta-\frac{1}{2}\gamma^2-c_2\gamma)}}{(\sqrt{c_1^2+2\lambda}-c_1)(\sqrt{c_2^2+2\theta}-c_2)(\frac{1}{2}\gamma^2+c_1\gamma-\lambda)}\,.
\end{eqnarray*}
\end{corollary}

\subsection*{Acknowledgments}
The author would like to express his sincere thanks to Professors Krzysztof Dębicki and Peng Liu for valuable comments and remarks and especially for pointing the problems with random time intervals.

\end{document}